\documentclass[12pt,a4paper]{article}

\usepackage[dvips]{graphics}
\usepackage[dvips]{epsfig}

\newtheorem{theorem}{Theorem}

\newtheorem{lemma}[theorem]{Lemma}
\newtheorem{proposition}[theorem]{Proposition}
\newtheorem{remark}{Remark}
\newenvironment{proof}[1][Proof]{\textbf{#1.} }{\ \rule{0.5em}{0.5em}}

\usepackage{latexsym}
\usepackage{amsfonts}
\usepackage{amssymb}
\usepackage{amsmath}
\usepackage[english]{babel}

\def\G{\Gamma}
\def\g{\gamma}

\begin{document}

\title{Nonorientable 3-manifolds admitting coloured  triangulations with at most 30 tetrahedra
\footnote{Work performed under the auspices of G.N.S.A.G.A. of
C.N.R. of Italy.}}

\author{Paola Bandieri,\ Paola Cristofori \ and Carlo Gagliardi}

\maketitle

\begin{abstract}
We present the census of all non-orientable, closed, connected
3-manifolds admitting a rigid crystallization with at most 30
vertices. In order to obtain the above result, we generate,
manipulate and compare, by suitable computer procedures, all rigid
non-bipartite crystallizations up to 30 vertices.

\end{abstract}

\section{Introduction}

Within the study of 3-dimensional manifolds, it is often useful to
have significative examples to formulate and test conjectures, to
obtain classification results or to investigate patterns for
3-manifolds.

During the last ten years, several papers have been published
containing tables (\textit{censuses}) of 3-manifolds, satisfying
certain conditions. The criterion, which is usually adopted, is to
bound the possible number of tetrahedra in a triangulation of the
manifold. First Matveev presented the census of closed orientable
irreducible 3-manifolds having a triangulation formed by at most
six tetrahedra (\cite{M1}). More precisely, Matveev's results are
based on the representation of 3-manifolds by special spines and
his bound is the complexity of the manifold, i.e. the minimal
number of vertices in a special spines of the manifold, which
coincides (excluding some very particular cases) with the number
of tetrahedra in a minimal triangulation.

The orientable censuses were later extended by Ovchinnikov up to
complexity 7 (\cite{O}), by Martelli and Petronio up to complexity
10 (\cite{MP1}, \cite{MP2} ) and by Matveev himself up to complexity
11 (\cite{M3}, \cite{M4} ).

With regard to the non-orientable case, the first tables were made
by Amendola and Martelli up to complexity 7 (\cite{AM1},
\cite{AM2}) and Burton up to 7 tetrahedra (\cite{B1}, \cite{B2});
recently Burton completed the census up to 10 tetrahedra
(\cite{B3}).

%However their approaches are different. Amendola and Martelli are
%more interested in the manifolds and their geometric structures,
%than in listing all possible triangulations of a given complexity.
%Moreover, their generation process, unlike the other authors,
%exploits theoretical constructions and makes no use of the
%computer.

%On the other hand,
%In particular, Burton focused on minimal triangulations and
%their combinatorial structures, and generated, by computer, all
%possible triangulations of closed, irreducible and $\mathbb
%P^2$-irreducible manifolds, having a fixed (not higher than ten)  %greater ??
%number of tetrahedra. Burton's program, \textit{Regina},
%can construct "`pieces"' of triangulations, which are topologically
%significative parts of known minimal triangulations; it also
%allows to manipulate triangulations and compare them with the
%already completed tables up to complexity ten.

In this paper, we share Burton's approach of listing and analyzing
all possible triangulations of closed $3$-manifolds, restricting
our considerations to "coloured triangulations" or, equivalently,
to their dual "edge-coloured graphs" (see \cite{FGG}, \cite{BCG}).

%Moreover, we make no restrictions with regard to irreducibility,
%due to the combinatorial objects we use. They are
%\textit{edge-coloured graph}, i.e. a combinatorial tool for
%representing PL-manifolds, which can be used in any dimension and
%without necessarily assuming connectedness and non-empty boundary,
%too (see \cite{FGG}, \cite{BCG}).

Edge-coloured graphs can be easily encoded by matrices and thus
manipulated by computer, in order to recognize topological
properties and compute invariants of the underlying manifolds or
in order to change triangulations by means of moves which preserve
the homeomorphism type of the represented manifolds.

Within the theory of edge-coloured graphs, several results have been
obtained in generating and classifying catalogues of closed
3-manifolds.

The adopted bounds are usually the number of vertices of the graph
(equivalently the number of tetrahedra of the coloured
triangulation) or the \textit{regular genus} of the graph, an
invariant whose minimal value coincides with the Heegaard genus of
the represented manifold.

In the first case, orientable catalogues were first produced and
analyzed by Lins up to 28 vertices (\cite{L2}; the classification
was completed in \cite{CA1}); later they were extended to 30 vertices in
\cite{CC}. Moreover, Casali in \cite{CA2} started the generation and study
of non-orientable catalogues, completing it up to 26 vertices.

In this paper, we extend the above result to the cases of
non-orientable closed manifolds representable by edge-coloured
graphs with 28 and 30 vertices. The generation procedure remains
as in \cite{CA2}, while the main point in the classification is an
algorithm, already introduced in \cite{CC}, for subdividing the
catalogues into classes so that the elements of each class
represent the same manifold.
%Moreover, as we will show in section
%..., the main result can be summarized in the following:
%\begin{theorem} ...\end{theorem}

As further step, we identified the manifolds represented by each
class, by computation of invariants, comparison with known
edge-coloured graphs and, in some cases, by constructing coloured
triangulations of manifolds in Burton's tables which matched our
representatives.

Finally, we found that there exist exactly thirty-three closed
non-orientable 3-manifolds, of which sixteen are prime, admitting
a coloured triangulation with at most 30 tetrahedra.  A precise
description of the above prime manifolds will be presented at the
end of section 4.

Existing catalogues of genus two orientable manifolds have been
generated and studied up to 34 (\cite{CA4},\cite{BGR}) and 42
tetrahedra (\cite{KN}). Presently, we are examining the
non-orientable case up to 42 tetrahedra: the related results will
be the subject of a forthcoming paper.

\bigskip

\section{Coloured triangulations of $3$-manifolds}

Throughout this paper, manifolds, when not otherwise specified, will always be closed and connected.

A coloured $n$-complex is a pseudocomplex (\cite {HW}) $K$ of
dimension $n$ with a labelling of its vertices by
$\Delta_n=\{0,\ldots,n\}$, which is injective on the vertex-set of
each simplex of $K$.

%Let us denote by $V(K)$ the vertex-set of $K$.

%If $\sigma^h$ is an $h$-simplex of $K$ then the \textit{disjoint
%star $std(\sigma^h,K)$} of $\sigma^h$ in $K$ is the pseudocomplex
%obtained by taking the disjoint union of the simplices of $K$
%containing $\sigma^h$ and identifying the $(n-1)$-simplices
%containing $\sigma^h$ together with all their faces.

%The \textit{disjoint link $lkd(\sigma^h,K)$} of $\sigma^h$ in $K$ is
%the subcomplex of $std(\sigma^h,K)$ formed by the simplices which
%don't intersect $\sigma^h.$

%From now on, we shall restrict our attention to the coloured
%complexes $K$, such that each $(n-1)$-simplex is a face of exactly
%two $n$-simplices and for each simplex $\sigma$ of $K$,
%$std(\sigma,K)$ is strongly connected.

\medskip

An \textit{$(n+1)$-coloured graph} is a pair $(\G,\gamma)$, where
$\G$ is a graph, regular of degree $n+1$, and $\gamma :
E(\G)\to\Delta_n$ a map which is injective on each pair of adjacent
edges of $\G$.

In the following, we shall often write $\G$ instead of
$(\G,\gamma)$.

For each $B\subseteq\{0,\ldots,n\}$, we call $B$\textit{-residues}
of $(\G,\gamma)$ the connected components of the coloured graph
$\G_B=(V(\G),\gamma^{-1}(B))$; given an integer $m\in\{1,\ldots,n\}$
we call $m$\textit{-residue} of $\G$ each $B$-residue of $\G$ with
$\# B=m$.

An isomorphism $\phi : \G\to\G'$ is called a \textit{coloured
isomorphism} between the $(n+1)$-coloured graphs $(\G,\gamma)$ and
$(\G',\gamma')$ if there exists a permutation $\varphi$ of
$\Delta_n$ such that $\varphi\circ\gamma=\gamma'\circ\phi$.

Coloured graphs are an useful tool for representing manifolds (see
\cite{BCG} and \cite{FGG} for a survey on this topic); in fact there
is a bijective correspondence between a particular class of $(n+1)$
- coloured graphs and the class of coloured triangulations of
$n$-manifolds.

A direct way to see this correspondence is to consider, for each
$(n+1)$-coloured graph $\G$, the coloured complex $K(\G)$ obtaining
by the following rule:
\begin{itemize}
\item [-] for each vertex $v$ of $\G$, take an $n$-simplex $\sigma (v)$ and label
its vertices by $\Delta_n$;

\item [-] if $v$ and $w$ are vertices of $\G$ joined by an $i$-coloured
edge ($i\in\Delta_n$), then identify the $(n-1)$-faces of $\sigma
(v)$ and $\sigma (w)$ opposite to the $i$-coloured vertex.

\end{itemize}

%Given a coloured complex $K$, a direct way to see this
%correspondence is to consider the coloured graph $\G$ which is the
%dual 1-skeleton of  $K=K(\G)$, i.e. the vertices of $\G$ are the
%barycenters of the $n$-simplices of $K(\G)$ and the edges of $\G$
%are the 1-cells dual of the $(n-1)$-simplices of $K(\G)$. Obviously
%the $(n-1)$-simplex dual to an $i$-coloured edge $e$ has its
%vertices labelled by $\hat{\imath}=\Delta_n-\{i\}$. In this way, a
%bijective correspondence is established between the $h$-simplices
%$(0\leqslant h\leqslant dim K(\G))$ of $K(\G)$ and the
%$(n-h)$-residues of $\G$; more precisely, if $\sigma^h$ is an
%$h$-simplex of $K(\G)$, whose vertices are labelled by $\{i_0,\ldots
%,i_h\}$, there is a unique $(n-h)$-residue $\Xi$ of $\G$ whose edges
%are coloured by $\Delta_n\smallsetminus\{i_0,\ldots ,i_h\}$ and such
%that $K(\Xi)=lkd(\sigma^h,K)$, where $K(\Xi)$ is the
%$(n-h-1)$-coloured complex dual to $\Xi$.

See \cite{FGG} for a more precise description of the involved
constructions.

If $M$ is a manifold of dimension $n$ and $\G$ an $(n+1)$-coloured
graph such that $|K(\G)|\cong M$, we say that $M$ is {\it
represented\/} by $\Gamma$.

If, for each $i\in\Delta_n$, $\Gamma_{\hat{\imath}}$ is connected
(equivalently the corresponding coloured triangulation $K(\Gamma)$
has exactly one $i$-coloured vertex for each $i\in\Delta_n$), then
both the $(n+1)$-coloured graph $\Gamma$ and the coloured
triangulation $K(\Gamma)$ are called \textit{contracted};
furthermore, if $\Gamma$ represents an $n$-manifold $M$, then it is
called a \textit{crystallization} of $M$. Note that $M$ is
orientable iff $\G$ is  bipartite.

%Several topological properties of $M$ correspond to combinatorial
%properties of each $(n+1)$-coloured graph $\G$ representing $M$; for
%example $M$ is orientable iff $\G$ is  bipartite.

We can construct a coloured graph representing the connected sum of
two $n$-manifolds $M^\prime$ and $M^{\prime\prime}$ starting from
their graphs. In fact let $\G^\prime$ and $\G^{\prime\prime}$ be
$(n+1)$-coloured graphs representing $M^\prime$ and
$M^{\prime\prime}$ respectively. Let $x$ be a vertex of $\G^\prime$
and $y$ a vertex of $\G^{\prime\prime}$, then the $(n+1)$-coloured
graph $\G=\Gamma^\prime\#\Gamma^{\prime\prime}$ obtained by removing
$x$ from $\G^\prime$ and $y$ from $\G^{\prime\prime}$ and by gluing
the "hanging" edges according to their colours, represents
$M^\prime\# M^{\prime\prime}$ (see \cite{FGG})\footnote{It is
well-known that, if both manifolds don't admit
orientation-preserving automorphisms, there exist two
non-homeomorphic connected sums. Each corresponds to requiring $x$
to belong to a fixed bipartition class in $V(\G^\prime)$ and
choosing $y$ in one of the two different bipartition classes of
$V(\G^{\prime\prime})$}.

%Recall that,
%if both $M$ and $M^\prime$ are orientable and both don't admit
%orientation-preserving automorphisms, there are two different
%connected sums.

\begin{remark}\label{somma}\emph{If $\G$ is a $(n+1)$-coloured graph representing a
$n$-manifold $M$ and if there are in $\G\ \ n+1$ edges
$\{e_0,\ldots,e_n\}$, one for each colour $i \in\Delta_n$, such
that $\G - \{e_0,\ldots,e_n\}$ splits into two connected
components, then it is easy to reverse the above procedure and
construct two $(n+1)$-coloured graphs $\G^\prime$ and
$\G^{\prime\prime}$, representing two $n$-manifolds $M^\prime$ and
$M^{\prime\prime}$ respectively, such that $\G = \G^\prime \#
\G^{\prime\prime}$, hence $M = M^\prime \#
M^{\prime\prime}$.}\end{remark}

%In the following, we will refer to the above condition as
%``condition (\#)"

\medskip

An important role, within the theory of coloured graphs, is played
by combinatorial moves (\textit{dipole moves}) which transform an
$(n+1)$-coloured graph representing an $n$-manifold into another
(usually non-colour isomorphic) $(n+1)$-coloured graph,
representing the same manifold.

%Given a vertex $v$ of a $(n+1)$-coloured graph $(\G ,\g)$ and a subset $\{i_1,\ldots,
%i_k\}\subset\Delta_n$, let us denote by $C_{i_1,\ldots, i_k}(v)$ the
%$\{i_1,\ldots, i_k\}$-residue of $(\G ,\g)$ containing $v$.

If $x,y$ are two  vertices of a $(n+1)$-coloured graph $(\G,\g )$
joined by $k$ edges $\{e_1,\ldots, e_k\}$ with $\g(e_h)=i_h$, for
$h = 1, \ldots, k$, then we call $\theta=\{x,y\}$ a
\textit{k-dipole} or a \textit{dipole of type k} in $\G$,
\textit{involving colours} $i_1,\ldots, i_k$, iff $x$ and $y$
belong to different $(\Delta_n - \{ i_1,\ldots, i_k\})$-residues
of $\G$.

%In this case $(\G^\prime,\g^\prime)$ is said to be obtained from $(\G,\g)$ by
%\textit{deleting the $k$-dipole} $\theta$ if it is obtained as
%%follows: \begin{itemize} \item [(a)] delete the vertices $x,y$ and
%all their incident edges. \item [(b)] If  $x_h$
%(resp. $y_h$) is the vertex $i_h$-adjacent to $x$ (resp. to $y$),
%$ h = 1, \ldots, k$, then join $x_h$ and $y_h$ by means of an
%$i_h$-coloured edge $f_h$.
%\end{itemize}

In this case a new $(n+1)$-coloured graph $(\G^\prime,\g^\prime )$ can be obtained from $\G$ by deleting $x,y$ and
all their incident edges and joining, for each $i\in\Delta_n-\{i_1,\ldots ,i_k\}$, the vertex $i$-adjacent to $x$ to the vertex $i$-adjacent to $y$; $(\G^\prime,\g^\prime )$ is said to be obtained from $(\G,\g )$ by \textit{deleting the $k$-dipole} $\theta$.
Conversely $(\G,\g )$ is said to be obtained from $(\G^\prime,\g^\prime)$ by \textit{adding the $k$-dipole}.

\medskip

By restricting ourselves to $3$-manifolds (in the following this will always be the case), we can introduce further moves.

Let $(\G,\g)$ be a $4$--coloured graph. Let $\Theta$ be a subgraph
of $\G$ formed by a $\{i,j\}$-coloured cycle $C$ of length $m+1$
and a $\{h,k\}$-coloured cycle $C^\prime$ of length $n+1$, having
only one common vertex $x_0$ and such that $\{i,j,h,k\} =
\{0,1,2,3\}$. Then $\Theta$ is called an \textit{(m,n)--dipole}.

If $x_1,\,x_m,\,y_1,\,y_n$ are the vertices respectively
$i,\,j,\,h,\,k$-adjacent to $x_0$, we define the $4$-coloured graph $(\G^\prime,\g^\prime)$
 \textit{obtained from} $\G$ \textit{by cancelling the
(m,n)--dipole}, in the following way:
\begin{itemize}
\item [1)] delete $\Theta$ from $\G$ and consider the
product $\Xi$ of the subgraphs $C-\{x_0\}$ and $C^\prime-\{x_0\}$;
\par\noindent
\item [2)] for each $s,s^\prime\in \{1, \dots, n\}$
(resp. for each $r,r^\prime\in \{1, \dots, m\}$), let $e$ be the
edge joining $y_s$ and $y_{s^\prime}$ (resp. $x_r$ and
$x_{r^\prime}$) in $\G$. If $\g (e) = c \in \{0,1,2,3\}$, then,
for each $t\in \{1, \dots, m\}$ (resp. for each $t\in \{1, \dots,
n\}$), join the vertices $(x_t, y_s)$ and $(x_t, y_{s^\prime})$
(resp. $(x_r, y_t)$ and $(x_{r^\prime}, y_t)$) by a $c$-coloured
edge in $\Xi$;
\par\noindent
\item [3)] for all $r\in\{1, \dots, m\},\,s\in\{1, \dots, n\}$, if a vertex $z$
of $\Gamma -\Theta$ is joined to $y_s$ (resp. $x_r$) by a $i$ or $j$ (resp.
$h$ or $k$)--coloured edge in $\G$,
then $z$ is joined to $(x_1,\,y_s),\,(x_m,\,y_s)$ (resp. $(x_r,\,y_1),\,
(x_r,\,y_n)$) by a $i$ or $j$ (resp. $h$ or $k$)--coloured edge in $\G'$.
\end{itemize}

\medskip

Figure 2.1 shows the whole process in the case $m=3$ and $n=5$.

\bigskip
\smallskip
\centerline{\scalebox{0.6}{\includegraphics{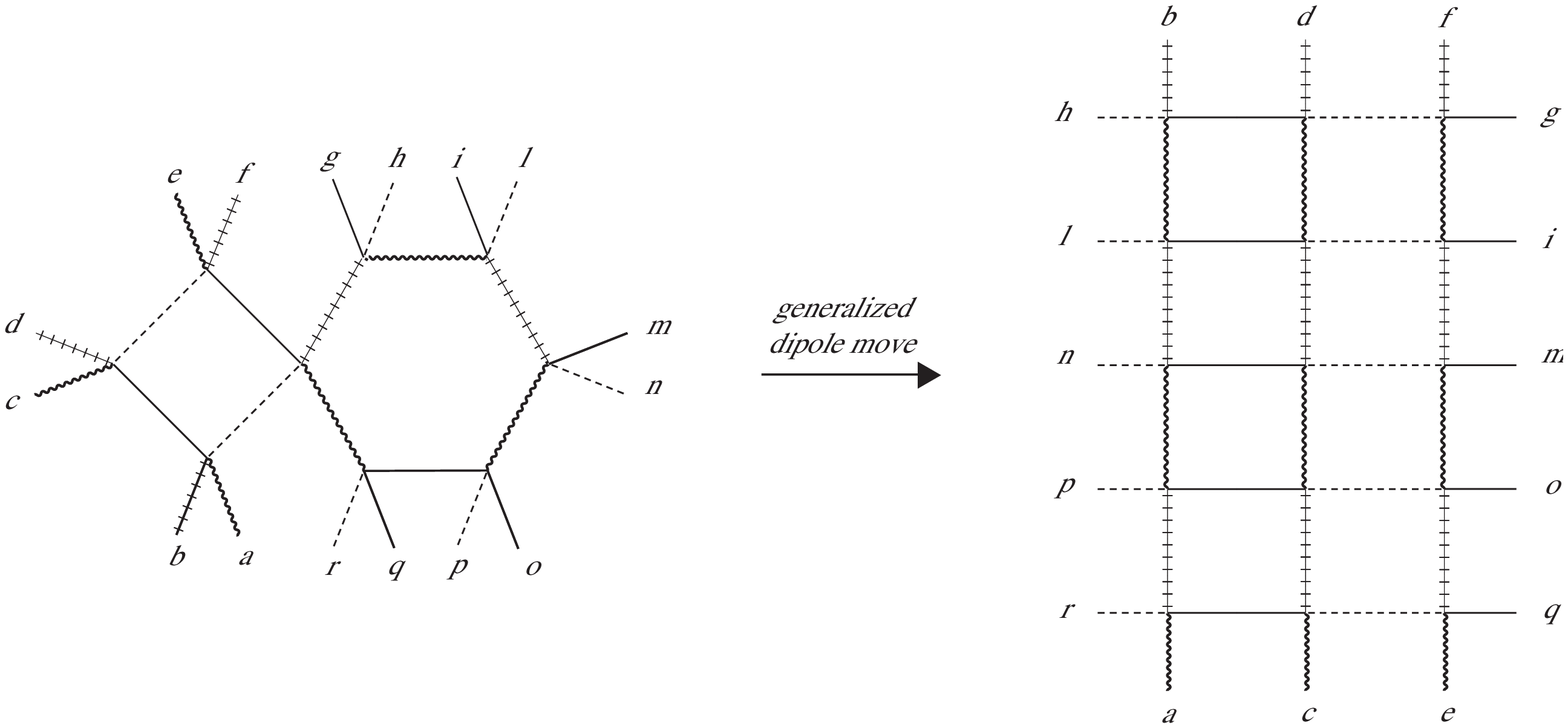}}}
\bigskip
\centerline{{\bf Figure 2.1}}

\bigskip

The moves just described are called \textit{generalized dipole
moves}. We can summarize in the following result the significance of dipole
moves and generalized dipole moves as a tool to manipule
4-coloured graphs.

\begin{proposition} {\rm (\cite{FG})} \ \label{dip-gen}
If $(\G,\g)$ and  $(\G^{\prime},\g^{\prime})$ are
4-coloured graph representing two 3-manifolds $M$ and $M^\prime $
respectively, and $( \G^{\prime},\g^{\prime})$ is obtained
from $(\G,\g)$ by a dipole move or a generalized dipole
move, then $M\cong M^\prime$.
\end{proposition}

Moreover there exist further moves, which can be applied to each
4-coloured graph $\G$, so to eliminate particular
configurations.

\medskip

If two $i$-coloured edges $e,f\in E(\Gamma)$ belong to the same
$\{i,j\}$-coloured cycle and to the same $\{i,k\}$-coloured cycle
of $\Gamma$, with  $j,k\in \Delta_3-\{i\}$ (resp. to the same
$\{i,h\}$-coloured cycle of $\Gamma$, for each $h\in
\Delta_3-\{i\}$), then $(e,f)$ is called a {\it $\rho_2$-pair}
(resp. a {\it $\rho_3$-pair}). Usually, we will write \textit{
$\rho$-pair} instead of $\rho_2$-pair or $\rho_3$-pair.

The graph $\G$ is a {\it rigid crystallization}
of a 3-manifold $M^3$ if it is a crystallization of $M^3$ and contains no
$\rho$-pairs.

A non-rigid crystallization $\G$ of a 3-manifold $M$ can always be
transformed into a rigid one by \textit{switching $\rho$-pairs}
(see
\cite{L2} % and Figure ....
) and cancelling the dipoles which could be created in the process.
The switching of a $\rho_2$-pair doesn't change the represented
manifold, while, for a $\rho_3$-pair, we have the following result.

\begin{lemma}\label{ro3coppie} {\rm (\cite{L2})} Let $\G$ be a 4-coloured
graph containing a $\rho_3$-pair, if $\G^\prime$, obtained from
$\G$ by switching it, is a crystallization of a 3-manifold $M$,
then $\G$ represents the 3-manifold $M\# H$, where $H=\mathbb
S^1\times\mathbb S^2$ iff $\G$ and $\G^\prime$ are both bipartite or both non
bipartite, otherwise $H=\mathbb S^1\tilde\times\mathbb S^2$.
\end{lemma}

Since each closed connected 3-manifold admits a rigid
crystallization (see \cite{CA1} for a detailed proof), we can
always require the rigidity condition to be satisfied with no loss
with regard to the represented manifolds.

\medskip

An essential tool to deal with coloured graphs by computer is
the \textit{code}, that is a numerical "string", which describes completely
the combinatorial structure of the coloured graph (see
\cite{CG} for definition and description of the related
\textit{rooted numbering algorithm}). More precisely, we can state

\begin{lemma} {\rm (\cite{CG})} Two $(n+1)$-coloured graphs are colour-isomorphic
iff they have the same code.\end{lemma}

\medskip

Therefore, by representing each coloured graph by its code, we can
easily reduced any catalogue of crystallizations to one containing
only non-colour-isomorphic graphs.

\medskip

A different description of triangulations for 3-manifolds,
including the coloured ones, can be found in \cite{B1}, \cite{B2}, by means of face pairings and gluing permutation selections;
moreover, each triangulation has its dual skeleton as an
associated $4$-valent graph. For a fixed triangulation of a
$3$-manifold the face pairing is intrinsic in the definition of
the triangulation and it is called the \textit{associated face
pairing}. In particular a coloured triangulation is a
triangulation in the sense of \cite{B1}, \cite{B2}, equipped
with the associated face pairing and having the identity map as
gluing permutation selection. The associated graph, if coloured by
associating to each edge the colour of the opposite vertex, is
exactly the $4$-coloured graph representing the coloured
triangulation.

\bigskip

\section{Generating and analysing catalogue $\tilde{\bf C}^{(30)}$}

\bigskip

In this section we will describe the generation and analysis of the
catalogue $\tilde {\bf C}^{(30)}$ of all non-isomorphic rigid
non-bipartite crystallizations with at most $30$ vertices. Since
the rigidity condition is not restrictive with regard to the
represented manifolds, the topological classification of the
crystallizations in $\tilde {\bf C}^{(30)}$ yields the list of all
non-orientable closed 3-manifolds admitting a coloured
triangulation with at most 30 tetrahedra; moreover the catalogue
$\tilde {\bf C}^{(30)}$ also yields the list of such
triangulations encoded into graphs.  %in ??

\medskip

In \cite{CA1} an algorithm is introduced, which for each $p
\in\mathbb N$, produces the archive $\mathcal C^{(2p)}$ (resp.
$\tilde {\mathcal C}^{(2p)}$) of codes of all non-isomorphic rigid
bipartite (resp. non-bipartite) crystallizations with exactly $2p$
vertices.

For our convenience, we summarize below the main steps of the algorithm.

\begin{itemize}
\item [Step 1: ] By induction on $p$ and by making use of the
results of \cite{L1} and \cite{L2}, we construct the set $\mathcal
S^{(2p)}= \{ \Sigma^{(2p)}_1, \Sigma^{(2p)}_2, \dots,
\Sigma^{(2p)}_{n_p} \}$ of all (connected) rigid 3-coloured graphs
with $2p$ vertices representing $\mathbb S^2$.

\item [Step 2: ] For each $i=1,2,\dots,n_p$, we add $p$ edges coloured by 3 to $\Sigma^{(2p)}_i$
in all ways so to obtain 4-coloured graphs, provided that the
planarity of the $\hat\imath$-residues ($i\in\{0,1,2\}$) and the
rigidity of the whole graph are preserved after each edge is
added.
Each time %Whenever  ??
a regular 4-coloured graph is obtained, we check whether it is a crystallization (i.e.
its Euler characteristic is zero).

\item[Step 3: ] By comparing the codes and by checking rigidity
condition and bipartition property on the crystallizations arising
from Steps 1 and 2, we form the catalogue ${\mathcal C}^{(2p)}$
(resp. $\tilde {\mathcal C}^{(2p)}$) of all rigid bipartite (resp.
non-bipartite) crystallizations with $2p$ vertices .

\end{itemize}

\bigskip

With regard to the non-bipartite case, the output data of a C++
program implementing the above algorithm, are shown in the
following Table (see \cite{CC} for the orientable case).

\bigskip

 \centerline{ \begin{tabular}{|c|c|c|c|c|c|c|c|c|c|c|c|c|c|c|c|}
  \hline
  %\ & \ & \ & \ & \ & \ & \ & \ & \ & \ & \ & \ & \ & \ & \ & \ & \\
  \hfill {\bf 2p } & 2 &  4 &  6 &  8 & 10 & 12 & 14 & 16 & 18 & 20 &
 22 & 24 & 26 & 28 & 30\\
 \hline
 %\ & \ & \ & \ & \ & \ & \ & \ & \ & \ & \ & \ & \ & \ & \ &\\
  {\bf $\#\tilde {\mathcal C}^{(2p)}$} & 0 & 0  & 0 & 0 & 0  & 0 & 1 & 1
   & 1 & 9 & 12 & 88 & 480 & 2790 & 21804\\
   %\ & \ & \ & \ & \ & \ & \ & \ & \ & \ & \ & \ & \ & \ & \ & \\
   \hline \end{tabular}}

\bigskip

\centerline{\textbf{Table 1:  non-bipartite rigid crystallizations
up to 30 vertices.}}

\bigskip

Catalogues $\tilde{\mathcal C}^{(2p)}$, for $p<14$, have been
analysed and the represented manifolds identified in \cite{CA1}
(see also \cite{CC1} and \cite{CC} for the orientable case),
mainly by manipulating the crystallizations through generalized
dipole moves and by subdividing them into classes according to the
equivalence defined by the moves.

In this paper we follow the same idea and generalize it into a more
refined ``classification" algorithm: more precisely, we will show
how to subdivide a given list $X$ of rigid crystallizations into
disjoint classes $\{c_1,\ldots ,c_s\}$ such that, for each
$i\in\{1,\ldots ,s\}$ and for each $\G ,\G^\prime\in c_i$,
there exist two integers $h,k\geq 0$ and a 3-manifold $M$ such that
$|K(\G)|=M\#_h H$ and $|K(\G^\prime)|=M\#_k H$, where $\#_r
H$ denotes the connected sum of $r$ copies either of the orientable
or of the non-orientable $\mathbb S^2$-bundle over $\mathbb S^1$;
more precisely $H=\mathbb S^1\times\mathbb S^2$ iff $\G$ and
$\G^\prime$ are both bipartite or both non-bipartite.

To make the algorithm clearer, let us introduce some definitions
and notations.

Let $\G$ be a rigid crystallization and suppose an ordering of
its vertices is fixed so that we can write
$V(\G )=\{v_1,\ldots,v_{2p}\}$; given an integer
$i\in\{1,2,3\}$, we denote by $\theta_i(\G )$ the rigid
crystallization obtained from $\G$ by subsequent cancellations
of $(m,n)$-dipoles of type $\{0,i\},$ according to the following
rules:
\begin{itemize}
\item [-] $m,n< 9$ (this condition is necessary in order to bound
the possible number of vertices of $\theta_i(\G)).$

\item [-] Generalized dipoles of type $\{0,i\}$ are looked for and cancelled
for increasing value of the integer $m\cdot n$ and by starting from
vertex $v_1$ up to $v_{2p}$, i.e if $\delta (v_i)$ is a $(m,n)-$generalized dipole at vertex $v_i$ (resp.
$\delta^\prime (v_j)$ is a $(m^\prime,n^\prime)-$ generalized dipole at vertex
$v_j$), then the cancellation of $\delta (v_i)$ is performed before
the cancellation of $\delta (v_j)$ iff ($m\cdot n<m^\prime\cdot
n^\prime$) or ($m\cdot n=m^\prime\cdot n^\prime$ and $i<j$).

\item [-] After each generalized dipole cancellation,
proper dipoles and $\rho$-pairs are cancelled in the resulting graph.
\end{itemize}

Moreover, we define $\theta_0(\G )=\G$.

\medskip

Given a rigid crystallization $\G$, there is a natural way to
construct a rigid crystallization $\G^{<}$ which is
colour-isomorphic to $\G$ and such that an ordering is induced
in $V(\G^{<})$ by the \textit{rooted numbering algorithm}
generating the code of $\G$ (see \cite{CG}). As a consequence,
for each $i\in\{0,1,2,3\}$, we can define a map $\theta_i$ on any
set $X$ of rigid crystallizations by setting, for each $\G\in
X$, $\theta_i(\G )=\theta_i(\G^{<})$, with the ordering of
the vertices induced by the code of $\G$.

\medskip

Let us denote by $S_3^0$ the set of all permutations on $\Delta_3$, which fix the element 0.
If  $S_3^0$ is considered as a lexicographically ordered set, let
$\delta^{(k)}=(\delta^{(k)}_0=0,
\delta^{(k)}_1,\delta^{(k)}_2,\delta^{(k)}_3)$ \
($k\in\{1,2,\ldots,6\}$) denote the $k$-th element of $S_3^0$ .

For each $k\in\{1,2,\ldots,6\}$ and for each $i\in\Delta_3,$ we
set
$$\ll \delta^{(k)}_i\gg=\theta_{\delta^{(k)}_i}\circ\theta_{\delta^{(k)}_{i-1}}\circ\ldots\circ\theta_{\delta^{(k)}_0}.$$

\medskip

Let us now consider the following set of moves:

$$\begin{aligned} \bar{\mathcal S} \, = & \,
\{\ll \delta^{(k)}_i\gg \ / \ k\in\{1,2,\ldots,6\} \text{ and }
i\in \Delta_3 \} \, \cup \\
\ & \cup \, \{\ll \delta^{(k)}_i\gg \circ \ll
 \delta^{(k-1)}_3\gg \circ \dots \circ \ll \delta^{(1)}_3\gg
\ / \ k\in\{2,\ldots,6\} \text{ and } i\in \Delta_3 \}
\end{aligned}, $$

\noindent and, for each rigid crystallization $\G$ and for each
$\epsilon\in\bar{\mathcal S}$, let $\theta_\epsilon (\G)$ be the
rigid crystallization obtained by applying the sequence of moves
$\epsilon$ to $\G$.

\medskip

Note that, by Propositions \ref{dip-gen} and \ref{ro3coppie}, each
sequence of moves $\epsilon\in\bar{\mathcal S}$ transforms a rigid
crystallization of a 3-manifold $M$ into a rigid crystallization
of a 3-manifold $M^\prime$, such that $M=M^\prime \#_t H$ ($H$ as
above) and $t$ is the number of $\rho_3$-pairs, which have been
deleted while performing the sequence $\epsilon$. As a
consequence, the following definitions naturally arise.

\medskip

For each $\G\in X$, the class $cl(\G )$ of $\G$ is defined as
\begin{align*}
cl(\G )=\{&\G^\prime\in X\ |\ \exists\,\epsilon,\epsilon^\prime\in
\bar{\mathcal S}\ \text{\ s.t.\ }\\
&\theta_\epsilon(\G ) \text{\ and\ }
\theta_{\epsilon^\prime}(\G^\prime) \text{\ have the same code}\}.
\end{align*}

Furthermore, we will denote by $h_\epsilon (\G)$ the number of
$\rho_3$-pairs which have been deleted by passing from $\G$ to
$\theta_\epsilon(\G)$ (obviously it could be zero).

Let us describe now the algorithm which, starting from an ordered
list $X$ of rigid crystallizations, simultaneously produces, for
each $\G\in X$, the set $cl(\G)$ and a non-negative number
$h(\G)$, whose meaning will be clear in the following.

More precisely, we will form $cl(\G )$ and compute $h(\G)$ in the following way.

\begin{itemize}
\item [Step 1: ] We set $cl(\G )=\{\G\}$ and
$h(\G )=0$;
\item [Step 2: ] for each $\epsilon\in\bar{\mathcal S}$, if there exist $\G^\prime\in X$
(coming before $\G$ in $X$) and $\epsilon^\prime\in\bar{\mathcal
S}$ such that the codes of $\theta_\epsilon(\G )$ and
$\theta_{\epsilon^\prime}(\G^\prime)$ coincide, then
\begin{itemize}
\item[$\bullet$] if $h(\G^\prime)-h_{\epsilon^\prime}(\G^\prime)\geq
h(\G )-h_{\epsilon}(\G )$, set $h(\G^{\prime\prime})=k-h(\G
)+h_{\epsilon}(\G )+h(\G^\prime)- h_{\epsilon^\prime}(\G^\prime)$
\ \ for each $\G^{\prime\prime}\in cl(\G )$ with
$h(\G^{\prime\prime})=k;$
\item[$\bullet$] if $h(\G^\prime)-h_{\epsilon^\prime}(\G^\prime)<
h(\G )-h_{\epsilon}(\G )$, set $h(\G^{\prime\prime})=k+h(\G
)-h_{\epsilon}(\G )-h(\G^\prime)+ h_{\epsilon^\prime}(\G^\prime)$
\ \ for each $\G^{\prime\prime}\in cl(\G^\prime)$ with
$h(\G^{\prime\prime})=k;$

\end{itemize}
\item [] In both cases, set $c=cl(\G)\cup cl(\G^\prime)$ and
$cl(\G^{\prime\prime})=c$, \ for each $\G^{\prime\prime}\in
c$.
\end{itemize}

\medskip
Furthermore, for each $c_i=\{\G_1^i,\ldots,\G_{r_i}^i\}$
and for each $0\leq h\leq max\{h(\G_1^i), \ldots,$
$h(\G_{r_i}^i)\}$, the class $c_i$ can be naturally subdivided
into subsets $c_{i,h}=\{\G^i_j\in c_i\ |\ h(\G^i_j)=h\}$.

\bigskip

By Propositions \ref{dip-gen} and \ref{ro3coppie}, it follows very easily that, if
$\G\in X$, $\G$ bipartite (resp. non-bipartite),
represents the manifold $M$ with
 $h(\G )=h$ and $c_i=cl(\G )$, then each
 element of $c_{i,k}$ ( $0\leq k\leq max\{h(\G^\prime)\ |\ \G^\prime\in c_i\}$)
  represents the manifold $M^\prime$ with $M^\prime=M\#_{k-h}H$
 or $M=M^\prime\#_{h-k} H$ according to $k \ge h$ or $k<h$ and $H$ as above.

\bigskip

\noindent\begin{remark}\emph{Note that the algorithm works as well
for any set $\mathcal S$ of sequences of generalized dipoles
moves, dipole moves and $\rho$-pairs switching, provided that each
element of $\mathcal S$ transforms rigid crystallizations into
rigid crystallizations; actually we could have described the above
procedure for such a general set $\mathcal S$ independently from
how the moves were performed (for this approach see \cite{CC1}).
However, we preferred to restrict ourselves to the particular set
$\bar{\mathcal S}$, which was used for our implementation.
}\end{remark}

\bigskip

It is clear that, if there exist $i\in\{1,\ldots ,s\}$ and $\G\in
c_i$ such that $|K(\G)|$ is known, then all manifolds represented
by crystallizations of $c_i$ are completely identified.
%; in fact supposing $\G\in c_{i,h}$ and $|K(\G)|=M\#_h H$ , then for each
%$0\leq k\leq max\{h(\G^\prime)\ |\ \G^\prime\in c_i\}$ and for
%each $\G^\prime\in c_{i,k}$, we have $|K(\G^\prime)|=M\#_k H$.
\medskip

Therefore, if known catalogues of crystallizations are inserted in
$X$, all classes of $X$ containing at least one known
crystallization are completely identified.

According to this idea the classification algorithm has been
implemented in the C++ program {\it $\G$-class}\footnote{{\it
$\G$-class} has been developed by M.R. Casali and P. Cristofori
and is available at WEB page
http://cdm.unimo.it/home/matematica/casali.mariarita/CATALOGUES.htm
where a detailed description of the program can be found, too}:
its input data are a list $X$ of rigid crystallizations and the
informations about the already known crystallizations of $X$
(possibly none), i.e the identification of their represented
manifolds through suitable ``names"; the output is the list of
classes of $X$, together with their representatives and, if
possible, their names.

We have applied {\it $\G$-class} to $X=\tilde{\textbf
C}^{(30)}=\bigcup_{1\leq p\leq 15} \tilde{\mathcal C^{(2p)}}$,
obtaining 32 classes: twelve contained crystallizations of
$\tilde{\mathcal C}^{(2p)}$ with $p<14$ and were therefore
completely identified. Furthermore four contained non-orientable
handles and were recognized by means of switching of
$\rho_3$-pairs and comparison with catalogues $\mathcal C^{(2p)}\
,\ 1\leq p\leq 15$ of rigid bipartite crystallizations.

Our further step was to apply Remark \ref{somma} to all
crystallizations of the unknown 16 classes, in order to recognize
possible connected sums: more precisely it was necessary to check
the condition of Remark \ref{somma} on each crystallization $\G$
of an unknown class and, in case of it being satisfied, to
construct the crystallizations $\G^{\prime}$ and
$\G^{\prime\prime}$ and try their recognition. This has also been
made by program {\it $\G$-class}, and the results involved 7
classes. They all represented distinct manifolds.

% except for nine classes all encoding $L(2,1)\#L(2,1)\#(\mathbb{RP}^2\times\mathbb S^1)$.

%This, as we will proof in the next section, was the only case when
%different classes turned out to represent the same manifold.

\bigskip

\section{Main results}

\bigskip

Before presenting our results with regard to the complete
identification of the manifolds, whose minimal (with regard to the
number of vertices) rigid crystallizations have 28 or 30 vertices,
let us introduce some notations.

We will denote by $TB(A)$ the torus bundle corresponding to the
matrix $A \in GL(2, \mathbb Z)$, i.e. the closed $3$-manifold
obtained as quotient of $T \times I$, by identifying the
bottom and top torus by the homeomorphism of $T^2$, induced by the
matrix $A$ (with respect to a fixed basis of $\pi_1(T^2)$). Recall
that $TB(A)$ is non-orientable if and only if $det A =-1$.

In \cite{CA2} a construction is described to obtain a $4$-coloured
graph representing $TB(A)$ for each $A\in GL(2, \mathbb Z)$.
Figures 4.1 and 4.2 sketch the whole process for the manifold
$TB\begin{pmatrix} 3 & 2\\2 & 1\end{pmatrix}$. In particular,
Figure 4.1 shows a subdivision of the boundary tori of $T\times I$
with the curves $c_i$ and $c_i^\prime$ ($i=0,1$), whose
identification produces the self-homeomorphism of $T$ defined by
the matrix $A$.

Figure 4.2 shows the boundary of a 3-dimensional simplicial
complex $K^\prime$, which is the first baricentric subdivision of
the cube $K$ representing $I\times I\times I$, triangulated by
subdividing the top and bottom squares as in Figure 4.1 and
performing the join on the boundary of $K$ from one vertex in its
interior.

$K^\prime$ is coloured by labelling $i$, for each $i=0,1,2,3$, the
vertices dual to the $i$-dimensional simplices of $K$. Finally the
coloured triangulation of $TB(A)$ is obtained by identifying, for
each $h=1,\ldots ,80$, the pairs of triangles on $\partial
K^\prime$, which in the Figure are labelled $h$ and $h^\prime$
respectively; the identification is obviously performed by
respecting the colouring of the vertices.

\vskip 30pt
\smallskip
\centerline{\scalebox{0.6}{\includegraphics{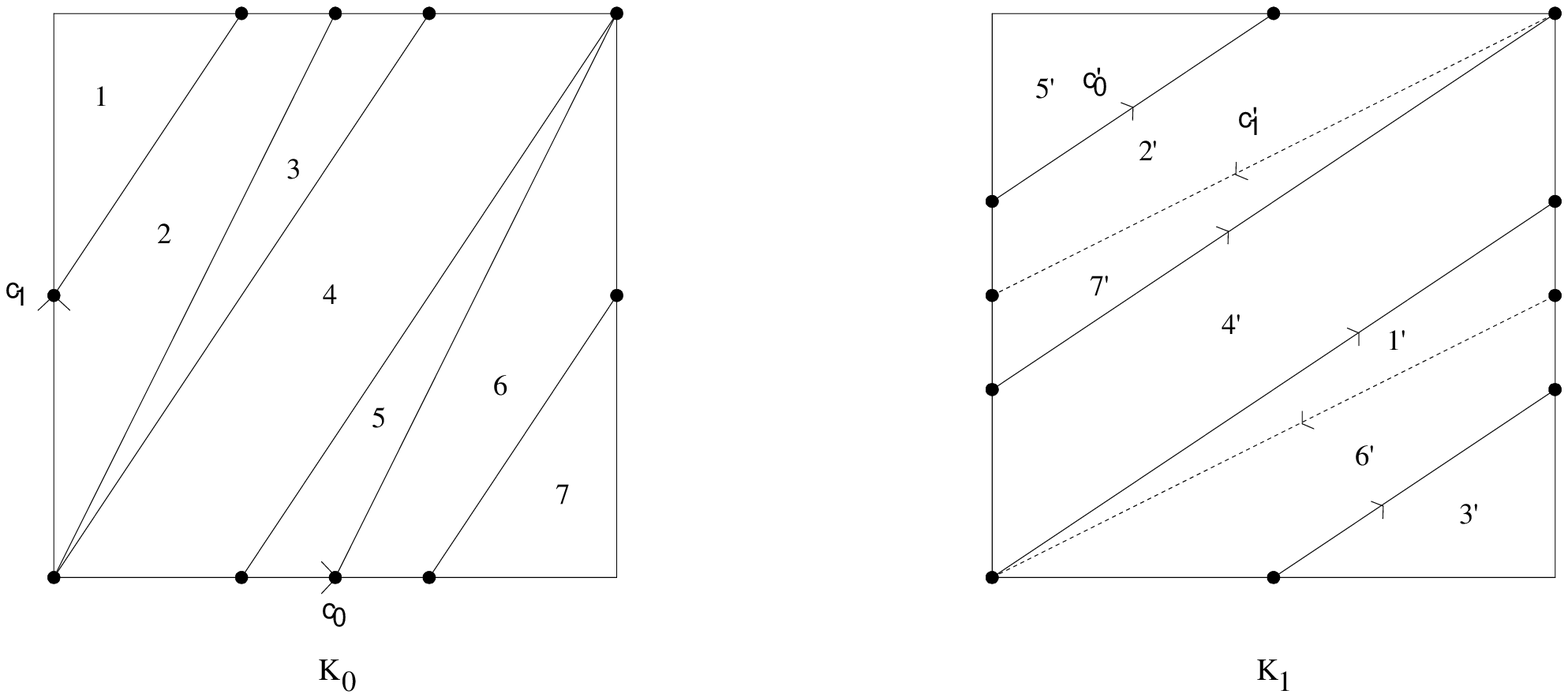}}}
\bigskip
\centerline{{\bf Figure 4.1}}

\bigskip

\bigskip
\centerline{\scalebox{0.9}{\includegraphics{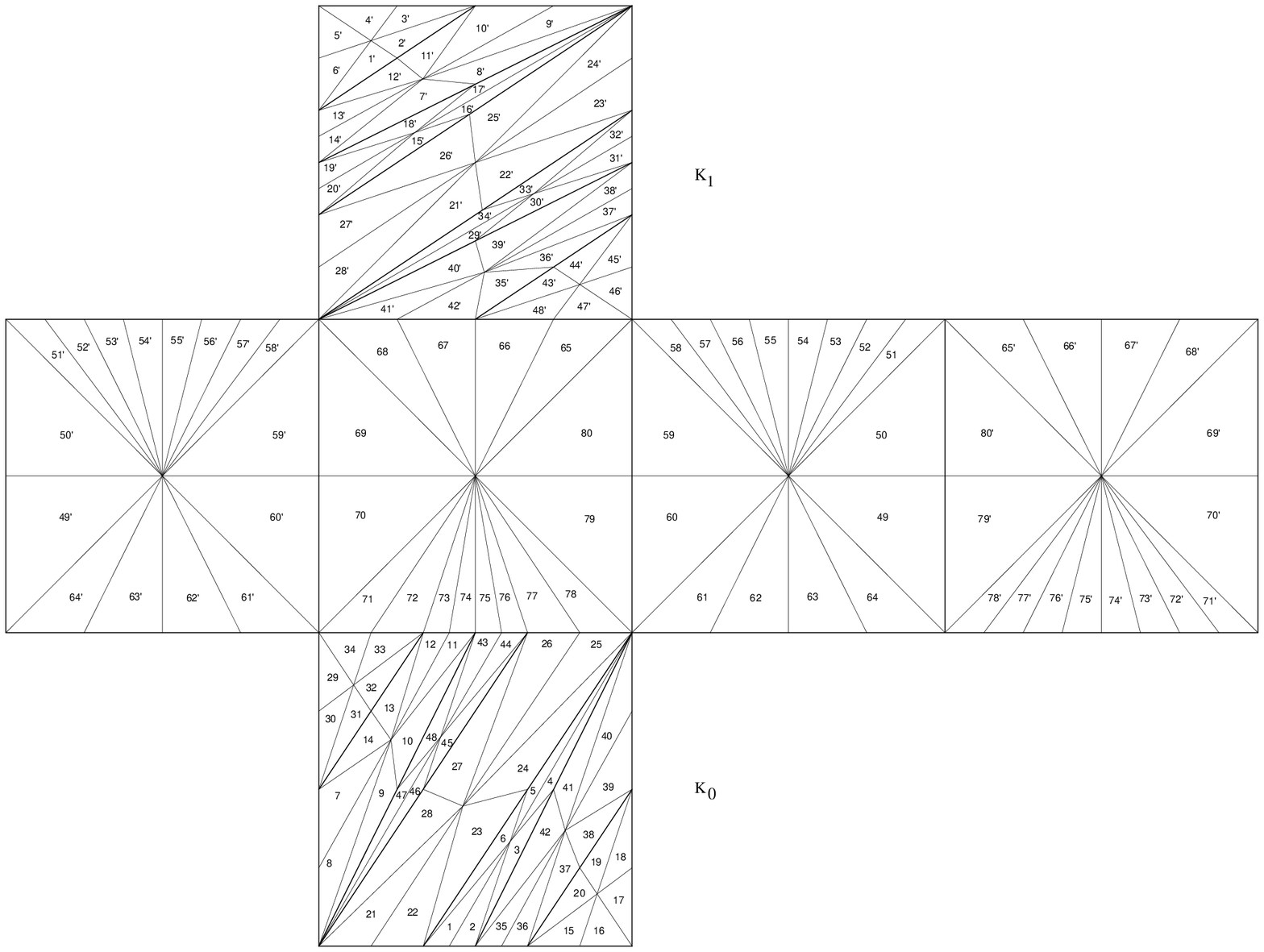}}}
\bigskip
\centerline{{\bf Figure 4.2}}

\vskip 20pt

Moreover, $SFS(S,(\alpha_1,\beta_1),\ldots ,(\alpha_r,\beta_r))$
will denote the Seifert fibred space with base space the orbifold
$S$ and (non-normalized) parameters $(\alpha_1,\beta_1),\ldots
,(\alpha_r,\beta_r)$ (in our case, we always have $r\leq 2$). The
base orbifolds of the Seifert manifolds in our catalogue, are

\begin {itemize}
\item [] $\mathbb{RP}^2\ $, projective plane
\item [] $\overline{D}\ $, disc with reflector boundary
\item [] $T^2/o_2\ $, torus containing fibre-reversing curves
\item [] $K^2\ $, Klein bottle
\item [] $K^2/n_3\ $, Klein bottle containing fibre-reversing curves with non-orientable total space
\item [] $A\ $, annulus.
\end{itemize}

Furthermore $SFS(A,(2,1))\cup SFS(A,(2,1))/N$ is the non-geometric graph-manifold obtained by pasting together two copies of $SFS(A,(2,1))$ along their boundary tori according to the matrix $N$ (which in our case is $\begin{pmatrix}0 & -1\\1 & 0\end{pmatrix}$).

\bigskip

By the results of program {\it $\G$-class}, we identified 23
classes: there remain nine unknown ones. Here is the list of codes
of their representatives, which are the first to appear in our
catalogue.

\medskip

\begin{itemize}
\item [$\G^{(1)}$]

\begin{footnotesize} CABFDEIGHLJKNMINDCMGFJLHEAKBJhKnHljbDgCfLdGEkiBMAeNacmFI\end{footnotesize}
\item [$\G^{(2)}$]
\begin{footnotesize} DABCHEFGKIJMLONKNFEDCBIHLOJGMAGliNOkADcofbKHjgLIhJManCEmBFd\end{footnotesize}
\item [$\G^{(3)}$]
\begin{footnotesize} DABCGEFJHIMKLONJNLEDCHGKOIFBMAMieKcJIobFDOAChmGnkNjBLgfdHaEl\end{footnotesize}
\item [$\G^{(4)}$]

\begin{footnotesize} CABFDEIGHLJKNMIMDCKGFJNHEBLAMIFBHjNlDfnhAkmdJiLcKebCGEag\end{footnotesize}
\item [$\G^{(5)}$]
\begin{footnotesize} EABCDIFGHLJKOMNLONGFEDCJIMAKHBMkIHNBlDCmbgjEGJfihnKodcAFOaeL\end{footnotesize}
\item [$\G^{(6)}$]
\begin{footnotesize} DABCGEFJHIMKLONJMOEDNHGKAIFBLCKehObIkcmEgDiCLGJnljMANfBaoFHd\end{footnotesize}
\item [$\G^{(7)}$]
\begin{footnotesize} DABCGEFJHIMKLONJMOEDNHGKAIFBLCKigOmIckbEhDeCLHFnljBNAfMaoJGd\end{footnotesize}
\item [$\G^{(8)}$]
\begin{footnotesize} DABCGEFJHIMKLONJMOEDNHGKAIFBLCKJgOmjcAMfHDeCLhFnlIBNkEbaoiGd\end{footnotesize}
\item [$\G^{(9)}$]
\begin{footnotesize} DABCGEFJHIMKLONJMOEDNHGKAIFBLCKFNiMOAndmGDjhBgoHlJbkCLEaIecf\end{footnotesize}
\end{itemize}

\medskip

The first step, in order to distinguish and recognize the manifolds
represented by $\G^{(1)},\ldots ,\G^{(9)}$, has been to write a
presentation of their fundamental groups by means of the algorithm
described in \cite{G}.

By abelianizing the presentation, we could compute the first
homology group of the involved manifolds. More precisely, we had

\medskip

$H_1(|K(\G^{(1)})|)=H_1(|K(\G^{(5)})|)=\mathbb Z\oplus
\mathbb Z_2$

$H_1(|K(\G^{(2)})|)=\mathbb Z\oplus \mathbb Z_3$

$H_1(|K(\G^{(3)})|)=\mathbb Z\oplus \mathbb Z_2\oplus \mathbb
Z_2$

$H_1(|K(\G^{(4)})|)=H_1(|K(\G^{(9)})|)=\mathbb Z$

$H_1(|K(\G^{(6)})|)=H_1(|K(\G^{(7)})|)=\mathbb Z\oplus
\mathbb Z$

$H_1(|K(\G^{(8)})|)=\mathbb Z\oplus \mathbb Z_8$.

\medskip

The above list allows us to establish that $\G^{(2)},\
\G^{(3)},\ \G^{(8)},\ $ represent distinct manifolds,
which are distinct from all others.

With regard to $\G^{(1)},\ \G^{(2)},\ \G^{(3)},\ $
however, a close analysis of the fundamental groups led us to
their identification as follows.

Let us consider the following presentation of
$\pi_1(|K(\G^{(1)})|)$, which comes from the algorithm in
\cite{G} by choosing the colours $\{1,3\}$:
$$<\ a,b,c,d,e\ /\ a^{-1}cd^{-1},\ a^{-1}be,\ c^{-1}ed^{-1}b, c^{-1}da,
e^{-1}ab^{-1}da^{-1}\ >$$

We get $d$ from the first and fourth relation, $e$ from the second relation
and substitute them in the remaining ones; hence we obtain the new
presentation

$$<\ a,b,c\ /\ [a^{-1},c]=1,\ bcb^{-1}ca^{-2},\
ba^{-1}b^{-1}a^{3}c^{-1}\ >$$

\noindent which is easily recognized as a presentation of $\mathbb
Z\begin{pmatrix} 3 & 1\\ -2 & -1\end{pmatrix}$, i.e. the
semidirect product of $\mathbb Z$ and $\mathbb Z\times\mathbb Z$
induced by the indicated matrix; hence
$\pi_1(|K(\G^{(1)})|)=\pi_1(TB\begin{pmatrix} 3 & 1\\ -2 &
-1\end{pmatrix})=\pi_1(TB\begin{pmatrix} 2 & 1\\ 1 &
0\end{pmatrix}$), since the last two matrices are conjugated in
$GL(2,\mathbb Z)$.

Analogously, if we consider the following presentation of
$\pi_1(|K(\G^{(2)})|)$ (the chosen colours are $\{0,2\}$):

$$<\ a,b,c,d,e\ /\ a^{-1}bd^{-1}ec,\ c^{-1}db,\ d^{-1}ae^{-1}, a^{-1}ec^{-1}d,
d^{-1}b^{-1}c\ >$$

\noindent by using the second, fifth and third relations to obtain
$d$ and $e$ and rewriting the remaining relations only in $a,\ b,\
c $, we have

$$<\ a,b,c\ /\ [b^{-1},c]=1,\ aca^{-1}b^3c^{-2},\
ab^{-1}a^{-1}bc^{-1}\ >$$

\noindent which is $\mathbb Z\begin{pmatrix} 1 & 1\\ 3 & 2\end{pmatrix}$.
Again, we have that $|K(\G^{(2)})|$ has the same fundamental
group as the torus bundle $TB\begin{pmatrix} 1 & 1\\ 3 &
2\end{pmatrix}=TB\begin{pmatrix} 3 & 1\\ 1 & 0\end{pmatrix}$.

Finally, we can perform similar transformations starting from

$$\pi_1(|K(\G^{(3)})|)=<\ a,b,c,d,e\ /\ b^{-1}ca^{-1}e^2,\
d^{-1}ae,\ d^{-1}ea, a^{-1}ec^{-1}eb, a^{-1}bc^{-1}\ >$$

\noindent and obtain

$\pi_1(|K(\G^{(3)})|)=\pi_1(TB\begin{pmatrix} 3 & 2\\ 2 &
1\end{pmatrix})$.

\medskip

The former analysis of fundamental groups led us to check the
suspected identifications; as a consequence we established the
following

\begin{proposition} $|K(\G^{(1)})|=TB\begin{pmatrix} 2 & 1\\ 1 &
0\end{pmatrix}$

$|K(\G^{(2)})|=TB\begin{pmatrix} 3 & 1\\ 1 & 0\end{pmatrix}$

$|K(\G^{(3)})|)=TB\begin{pmatrix} 3 & 2\\ 2 & 1\end{pmatrix}$

\end{proposition}

\begin{proof} By means of the algorithm described in \cite{CA2},
we constructed triangulations of $\ TB\begin{pmatrix} 2 & 1\\ 1 &
0\end{pmatrix},\ TB\begin{pmatrix} 3 & 1\\ 1 & 0\end{pmatrix},\
TB\begin{pmatrix} 3 & 2\\ 2 & 1\end{pmatrix}$ (see Figure 4.3 for
the last) and, by means of dipoles cancellation and switching of
$\rho$-pairs\footnote{This process was performed by using program
Duke III, available at WEB page
http://cdm.unimo.it/home/matematica/casali.mariarita/DukeIII.htm,
which manipulates edge-coloured graphs.}, we obtained rigid
crystallizations of the above torus bundles, which were added to
the set $cl(\G^{(1)})\cup cl(\G^{(2)})\cup
cl(\G^{(3)})$. The classification program applied to this list
produced exactly three classes, proving the above
identifications.\end{proof}

\bigskip

In order to analyze the remaining six unknown classes, we used a
combinatorial invariant of 3-manifold, the \textit{GM-complexity}
which is an upper bound for Matveev's complexity (\cite{M1},
\cite{M2}).

The GM-complexity $c_{GM}(M)$ of a 3-manifold $M$ is defined as $\
\min c(\G)\ $, where the minimum is taken among all
crystallizations $\G$ of $M$ and $c(\G)$, the
\textit{complexity} of the 4-coloured graph $\G$, is an
integer which can be computed on $\G$ by means of an easily
implemented algorithm (see \cite{CA3} and \cite{CC1} for precise
definition and results).

By the computation of complexity performed on all representatives
of $cl(\G^{(4)}),\ldots ,cl(\G^{(9)})$\footnote{Again GM-complexity
computation was performed by a C++ program, {\it CGM} (authors M.R. Casali and
P. Cristofori), available at WEB page
http://cdm.unimo.it/home/matematica/casali.mariarita/DukeIII.htm}, we had

\medskip

$c_{GM}(|K(\G^{(i)}|)\leq 8,\ \ $ for $i=4,5$

$c_{GM}(|K(\G^{(i)}|)\leq 9,\ \ $ for $i=6,7,8,9$.

\medskip

Table 9 of \cite{B3}, which shows the closed non-orientable
3-manifolds up to Matveev's complexity 10, gave us
possible identifications for $|K(\G^{(i)})|,\ \ i=4,\ldots,9$
and the program \textit{Regina}, realized by B. Burton
\footnote{the program \textit{Regina} is available at WEB page
http://regina.sourceforge.net} gave us the list of minimal
triangulations of these manifolds.

\medskip

It is easy to see that, given a triangulation $T$ of a manifold
$M$, a coloured triangulation of $M$ can be always constructed by
taking the first barycentric subdivision $T^\prime$ of $T$ and
labelling each vertex $v\in V(T^\prime)$ by the dimension of its
dual simplex of $T$.

Moreover the 4-coloured graph dual to this coloured triangulation
can be always reduced, by deleting dipoles and switching of
$\rho$-pairs, to a rigid crystallization of $M$.

Therefore, Burton's triangulations allowed us to obtain a list $X$
of cystallizations of manifolds, which were possible candidates for
our identifications: the set $\bigcup_{i=4}^9 cl(\G^{(i)})\cup X$
was handled by {\it $\G$-class}, whose results are summarized in
the following proposition.

\begin{proposition} The classes $cl(\G^{(4)},\ldots,cl(\G^{(9)})$
all represent distinct manifolds, which are precisely:

\medskip

$|K(\G^{(4)})|=SFS(\mathbb{RP}^2,(2,1),(3,1))$

$|K(\G^{(5)})|=SFS(\overline D,(2,1),(3,1))$

$|K(\G^{(6)})|=SFS(T^2/o_2,(2,1))$

$|K(\G^{(7)})|=SFS(K^2,(2,1))$

$|K(\G^{(8)})|=SFS(K^2/n_3,(2,1))$

$|K(\G^{(9)})|=SFS(A,(2,1))\cup SFS(A,(2,1))/\begin{pmatrix} 0
& -1\\ 1 & 0\end{pmatrix} $.
\end{proposition}

As a consequence of the above proposition and the results of the
previous section, we can state

\begin{proposition} There exists a one-to-one correspondence between the
set of classes of $\tilde{\textbf C}^{(30)}$ produced by the
classification program and the set of non-orientable 3-manifolds
admitting a coloured triangulation with at most 30 tethraedra.
\end{proposition}

\medskip

The following Table summarizes our results with regard to prime
non-orientable manifolds admitting a coloured triangulation with
at most 30 tethraedra \footnote{We recall that
%$\mathbb S^1\\widetilde{\times}\ \mathbb S^2$ is the twisted $\mathbb S^2$-bundle over $\mathbb S^1$, while
$\mathbb E^3/\mathcal B_i$ with $i=1,2,3,4\ $ denote the four
non-orientable Euclidean manifolds according to the notations of
\cite{V}.}.

%\smallskip

%\begin {theorem}\label{Main th} There are exactly sixteen prime, closed, nonorientable
%$3$-manifolds admitting coloured triangulations up to 30 tetrahedra; they are collected %resumed  ??
%in the following Table: \end{theorem}

\bigskip

 \centerline{ \begin{tabular}{|c|c|}
  \hline
  \hfill  \textbf{tethraedra} \hfill &
\textbf{$3$-manifold }  \\
 \hline \ & \ \\
  14 & $\mathbb S^1\ \widetilde{\times}\ \mathbb S^2$  \\ \ & \ \\\hline \ & \ \\ 16 & $\mathbb R\mathbb P^2 \times \mathbb S^1$   \\ \ & \ \\\hline \ & \ \\
  24 & $\mathbb E^3 /\mathcal B_1$ \\   & $\mathbb E^3 /\mathcal B_2$
 \\   & $\mathbb E^3 /\mathcal B_4$  \\ \ & \ \\\hline\ & \ \\ 26 & $\mathbb E^3 /\mathcal B_3$  \\   & $TB\begin{pmatrix} 0 & 1\\ 1 & -1\end{pmatrix}$ \\ \ & \ \\\hline \ & \ \\ 28 & $TB\begin{pmatrix} 2 & 1\\ 1 &
0\end{pmatrix}$   \\   & SFS($\mathbb R \mathbb P^2;
 (2,1),(3,1))$
 \\ \ & \ \\\hline \ & \ \\ 30 & $TB\begin{pmatrix} 3 & 2\\ 2 & 1\end{pmatrix}$\\  & SFS($\overline{D};(2,1),(3,1))$\\  & $TB\begin{pmatrix} 3 & 1\\ 1 & 0\end{pmatrix}$\\  & SFS($T^2/o_2;(2,1))$\\ & SFS($K^2;(2,1))$\\  &
 SFS($K^2/n_3;(2,1))$\\ & SFS($A;(2,1))\cup SFS(A;(2,1))/\begin{pmatrix}0 & -1\\1 & 0\end{pmatrix}$
\\ \ & \ \\
\hline \end{tabular}}

\bigskip
\centerline{\bf TABLE 2:  Prime 3-manifolds represented by
crystallizations of $\tilde{\textbf C}^{(30)}$}

\end{document}